\documentclass[12pt]{article}
\usepackage{amsmath}
\usepackage{amscd}
\usepackage{mathrsfs}
\usepackage{amsfonts}
\usepackage{amsmath,amssymb}
\baselineskip 5.5pt \lineskip 5.5pt \numberwithin{equation}{section}

\def\ad{\mbox{ad}}

\def \C{\hbox{$C\hskip -5pt \vrule height 6pt depth 0pt \hskip 6pt$}}

\def\qed{\ \ \ifhmode\unskip\nobreak\fi\ifmmode\ifinner
         \else\hskip5pt\fi\fi
 \hbox{\hskip5pt\vrule width4pt height6pt depth1.5pt\hskip 1 pt}}
\def\a{\alpha}
\def\Span{{\rm span}}

\def\d{\delta}
\def\D{\Delta}

\def\gi{\mathfrak{g}}

\def\l{\lambda}

\def\cl{\centerline}

\def\Der{\mathrm{Der}}

\def\D{\Delta}

\def\vs{\vspace*}

\def\C{\mathbb{C}}

\newtheorem{theo}{Theorem}[section]
\newtheorem{lemm}[theo]{Lemma}

\begin{document}
\cl {{\large\bf
 \vs{10pt} $2$-local derivation on the conformal Galilei algebra}
\noindent\footnote{Supported by the National Science Foundation of
China (Nos. 11047030 and 11771122).}} \vs{6pt}

\cl{Yufang Zhao, Yongsheng Cheng}
\cl{ \small School
of Mathematics and Statistics, Henan
University, Kaifeng 475004, China}

\vs{6pt}

{\small
\parskip .005 truein
\baselineskip 10pt \lineskip 10pt

\noindent{{\bf Abstract}\,  2-local derivation is a generalized derivation for a Lie algebra,
which plays an important role to the study of local properties of the structure of the Lie algebra.
In this paper, we prove that every 2-local derivation on the conformal Galilei algebra
is a derivation.
 \vs{5pt}

\noindent{\bf Key words}\,} derivation, 2-local derivation, the conformal Galilei algebra

\noindent{\bf MR(2000) Subject Classification} 16E40, 17B56, 17B68.
\parskip .001 truein\baselineskip 8pt \lineskip 8pt

\vs{6pt}
\par
\cl{\bf\S1. \ Introduction}
\setcounter{section}{1}\setcounter{theo}{0}\setcounter{equation}{0}

Let $g$ be an algebra and $\D$ a map of $g$ into itself. In \cite{[S]}, P. \v{S}emrl defined
$\D$ (not necessarily linear) to be a $2$-local derivation of $g$, if for every pair
of elements $x,y\in g$, there exists a derivations
$\D_{x,y}: g\rightarrow g$ (depending on $x, y$) such that
$\D_{x, y}(x)=\D(x)$ and $\D_{x, y}(y)=\D(y)$.
Furthermore, P. \v{S}emrl proved that every 2-local derivation of $B(H)$ is a derivation of $B(H)$,
where $H$ is an infinite-dimensional separable Hilbert space  and
$B(H)$ is the algebra of all bounded linear operators on $H$.
Clearly, as a generalization of derivation, \v{S}emrl introduced the notion of 2-local derivations on algebras.
The main problems concerning the notion are to find conditions under which 2-local derivations become derivations
and to present examples of algebras with 2-local derivations that are not derivations.

Investigation of 2-local derivations on finite dimensional Lie algebras and infinite
dimensional Lie (super) algebras were initiated in papers \cite{[AK1], [AK2], [DGL], [T], [ZCZ]}.
In \cite{[AK1]}, the authors proved that every 2-local derivation on a semi-simple Lie algebra is a derivation and that each finite-dimensional
nilpotent Lie algebra with dimension larger than two admits 2-local derivation which is not a derivation. In \cite{[AK2], [DGL], [T], [ZCZ]},
the authors proved that 2-local derivations on the Witt algebra, super Virasoro algebra, W-algebra $W(2, 2)$ and its superalgebra are derivations
and there are 2-local derivations on the so-called thin Lie algebra which are not derivations.

In this paper, we will study 2-local derivations on the conformal Galilei algebra.
The plan of this paper is as follows. In section 2, we give some preliminaries concerning the conformal
Galilei algebra, and determine the general form of the derivations on the conformal Galilei algebra. In section 3, we prove that
every 2-local derivation on the conformal Galilei algebra automatically becomes a derivation.

\cl{\bf\S2.\ Notations and Preliminaries}
\setcounter{section}{2}\setcounter{theo}{0}\setcounter{equation}{0}
In this paper, we denoted
 by $\mathbb{Z},\mathbb{N}$,and $\mathbb{C}$ the sets of all integers,
 positive integers, and complex numbers, respectively.
 In this section we give some necessary notations, definitions and preliminary results.

A derivation on a Lie algebra $g$ is a linear map $D:\, g\rightarrow g$ satisfying
$$D[x, y]=[D(x), y]+[x, D(y)]$$
for all $x, y\in g$. Denote by $Der(g)$ the set of all derivations of $g$.
For all $a\in g$, the map $ad(a)$ on $g$ defined as $ad(a)(x)=[a,x]$, $x\in g$ is a derivation
and derivations of this form are called inner derivation and denote by $ad(g)$.

A 2-local derivation is a natural generalization of a derivation of $g$. Clearly, for a 2-local derivation
on $g$ and $k\in\C$, $x\in g$, we have
$$\D(kx)=\D_{x,kx}(kx)=k\D_{x,kx}(x)=k\D(x).$$

In recent years Galilei groups and their Lie algebras have been intensively studied. The
interest is due to an appearance of this kind of symmetries in very different areas of physics and
mathematics \cite{[AGM], [GM]}. The conformal extension of the Galilei algebra is parameterized
by a positive half-integer number $l$ and is called the $l$-conformal~Galilei~algebra.
For
$l\in\mathbb{N}-\frac{1}{2}$,~we denote the conformal Galilei algebra by $\mathfrak{g}$, which has a basis given by
$$\{e,h,f,p_{k},z\mid k=0,1,2,\ldots,2l\}$$
and the Lie bracket given by:
$$\begin{array}{cc}
[h,e]\!=\!\!2e,~~~~~~[h,f]\!=\!\!-2f,~~~~~~[e,f]\!=\!\!h,\\[6pt]
     [h,p_{k}]\!=\!\!2(l-k)p_{k},~~~~~~[e,p_{k}]\!=\!\!kp_{k-1},
     ~~~~~~[f,p_{k}]\!=\!\!(2l-k)p_{k+1},\\[6pt]
     [z,\mathfrak{g}]\!=\!\!0,\\[6pt]
    [p_{k},p_{k'}]\!=\!\!\delta_{k+k',2l}(-1)^{k+l+\frac{1}{2}}k!(2l-k)!z, ~~\mathrm{for}~k,k'=0,1,2,\ldots,2l.
\end{array}$$
The conformal Galilei algebra can be viewed as a semidirect product
$\mathfrak{g}=sl_{2}\ltimes H_{l}$ of two subalgebras.  $sl_{2}=\Span\{e,h,f\}$ and Heisenberg
subalgebra $H_{l}=\Span\{p_{k},z\mid
k=0,1,2,\ldots,2l\}$. The irreducible representations of the conformal Galilei algebra are classified in \cite{[AI], [CCS], [LMZ], [ZC]}.

The following lemma can be used to determine the derivation of $\mathfrak{g}$.
\begin{lemm}\label{lemm0201}
Let $D$ be a linear map from $\mathfrak{g}$ into itself. Then $D\in Der (\mathfrak{g})$ if and
only if for $k,k'\in 0,1,\cdots,2l$ the following
equations hold:
\begin{eqnarray}\label{der-0}[D(h),e]+[h,D(e)]\!\!\!&=&\!\!\!2D(e) , \\[6pt]
\label{der-1}[D(h),f]+[h,D(f)]\!\!\!&=&\!\!\!-2D(f), \\[6pt]
\label{der-2}[D(e),f]+[e,D(f)]\!\!\!&=&\!\!\!D(h),\\[6pt]
\label{der-3}[D(h),p_{k}]+[h,D(p_{k})]\!\!\!&=&\!\!\!2(l-k)D(p_{k}),  \\[6pt]
\label{der-4}[D(e),p_{k}]+[e,D(p_{k})]\!\!\!&=&\!\!\!kD(p_{k-1}), \\[6pt]
\label{der-5}[D(f),p_{k}]+[f,D(p_{k})]\!\!\!&=&\!\!\!(2l-k)D(p_{k+1}),\\[6pt]
\label{der-6}[D(p_{k}),p_{k'}]+[p_{k},D(p_{k'})]\!\!\!&
=&\!\!\!\delta_{k+k',2l}(-1)^{k+l+\frac{1}{2}}k!(2l-k)!D(z).
\end{eqnarray}
\end{lemm}

Let $\delta$ be an outer derivation of $\mathfrak{g}$ determined by
\begin{equation}\label{out-01}
\delta(h)=\delta(e)=\delta(f)=0,~~~~~~\delta(z)=z,
~~~~~~\delta(p_{k})=\frac{1}{2}p_{k}, k=0, 1, \ldots, 2l.
\end{equation}
\begin{lemm}\label{daozi}
$\Der\gi=\C \d \oplus \mathrm{ad}(\gi)$, where $\d(p_i)=\frac{1}{2}p_i$,
 $\d(z)=z$ for $i=0,1,\dots,2l$.
\end{lemm}
\noindent{\it Proof.~}
Assume that $D\in Der \mathfrak{L}$ and $A=(a_{ij})_{(2l+5)\times(2l+5)}$
is the matrix of $D$ under the basis \{$e,h,f,p_{k},z\mid
k=0,1,2,\ldots,2l$\}, i.e.
\begin{equation}\begin{split}\label{der-7}
  &(D(p_{2l}), \cdots, D(p_{l+\frac{3}{2}}), D(f), D(p_{l+\frac{1}{2}}),
D(h), D(z), D(p_{l-\frac{1}{2}}), D(e), D(p_{l-\frac{3}{2}}), \cdots, D(p_{0}))\\
 =&(p_{2l},\cdots,p_{l+\frac{3}{2}},f,p_{l+\frac{1}{2}},h,z,p_{l-\frac{1}{2}},
 e,p_{l-\frac{3}{2}},\cdots,p_{0})A
\end{split}
\end{equation}
By (\ref{der-0}) and (\ref{der-7}), we have
\begin{multline}\label{der-8}
~~~~~~~a_{1,l+\frac{11}{2}}=a_{l+\frac{1}{2},l+\frac{11}{2}}=a_{l+\frac{7}{2},l+\frac{11}{2}}=a_{l+\frac{5}{2},l+\frac{5}{2}}=0,
~~~~~~~~(2-2(1-l))a_{2,l+\frac{11}{2}}=-2la_{1,l+\frac{5}{2}},\\
(2-2(l-k))a_{2l+1-k,l+\frac{11}{2}}=-(k+1)a_{2l-k,l+\frac{5}{2}}~~~~\mathrm{for}~~l+\frac{3}{2}<k<2l-1,\\
5a_{l-\frac{1}{2},l+\frac{11}{2}}=-(l+\frac{5}{2})a_{l-\frac{3}{2},l+\frac{5}{2}},
~~~~~~~~
3a_{l+\frac{3}{2},l+\frac{11}{2}}=-(l+\frac{3}{2})a_{l-\frac{1}{2},l+\frac{5}{2}},\\
2a_{l+\frac{5}{2},l+\frac{11}{2}}=-a_{l+\frac{1}{2},l+\frac{5}{2}},~~~~~~~~
a_{l+\frac{9}{2},l+\frac{11}{2}}=-(l+\frac{1}{2})a_{l+\frac{3}{2},l+\frac{5}{2}},\\
a_{l+\frac{13}{2},l+\frac{11}{2}}=(l-\frac{1}{2})a_{l+\frac{9}{2},l+\frac{5}{2}},
~~~~~~~~(2-2l)a_{2l+5,l+\frac{11}{2}}=-a_{2l+4,l+\frac{5}{2}},\\
(2-2(l-k'))a_{2l+5-k',l+\frac{11}{2}}=-(k'+1)a_{2l+4-k',l+\frac{5}{2}}~~~~\mathrm{for}~0<k'<l-\frac{3}{2}.~~~~
\end{multline}

Similarly, for $k\in 0,1,\cdots,2l$, ~~using (\ref{der-3})-(\ref{der-6}), (\ref{der-7}) and (\ref{der-8}), we obtain
 \begin{eqnarray*}
  \label{der-47}&&D(p_{2l})=\mu_{2}p_{2l}-2l\mu_{5}p_{2l-1}
+(-1)^{l+\frac{3}{2}}(2l-1)!\mu_{8}z\\[6pt]
  \label{der-48}&&D(p_{k})=(2l-k)\mu_{4}p_{k+1}+((2l-k)\mu_{1}+\mu_{2})p_{k}
-k\mu_{5}p_{k-1}+\frac{(-1)^{3l-k+\frac{1}{2}}(2l-k)!k!}{2(l-k)}\mu_{j_{2l-k}}z,\\[6pt]
  \label{der-49}&&D(p_{l+\frac{3}{2}})=(l-\frac{3}{2})\mu_{4}p_{l+\frac{5}{2}}
+((l-\frac{3}{2})\mu_{1}+\mu_{2})P_{l+\frac{3}{2}}-(l+\frac{3}{2})\mu_{5}p_{l+\frac{1}{2}}
-\frac{(l+\frac{3}{2})!(l-\frac{3}{2})!}{3}\mu_{10}z,\\[6pt]
  \label{der-50}&&D(f)=\frac{-\mu_{i_{2l-1}}}{-2+2l}p_{2l}+\cdots+
\frac{2l-k}{2(l-k)}\mu_{i_{k}}p_{k+1}+\cdots+(-\frac{l-\frac{3}{2}}{3}\mu_{9})p_{l+\frac{5}{2}}
+(l-\frac{1}{2})\mu_{6}p_{l+\frac{3}{2}}-\mu_{1}f\nonumber\\
  &&~~~~~~~~+(l+\frac{1}{2})\mu_{7}p_{l+\frac{1}{2}}-\mu_{5}h+\frac{l+\frac{3}{2}}{3}\mu_{10}p_{l-\frac{1}{2}}
+\frac{l+\frac{5}{2}}{5}\mu_{j_{l-\frac{5}{2}}}p_{l-\frac{3}{2}}+\cdots
+\frac{2l-k'}{2(l-k')}\mu_{j_{k'}}p_{k'+1}\nonumber\\
  &&~~~~~~~~+\cdots+\mu_{8}p_{1},\\[6pt]
  \label{der-51}&&D(p_{l+\frac{1}{2}})=(l-\frac{1}{2})\mu_{4}p_{l+\frac{3}{2}}
+((l-\frac{1}{2})\mu_{1}+\mu_{2})p_{l+\frac{1}{2}}-(l+\frac{1}{2})\mu_{5}p_{l-\frac{1}{2}}
+(l+\frac{1}{2})!(l-\frac{1}{2})!\mu_{7}z,\\[6pt]
  \label{der-52}&&D(h)=(-\mu_{3})p_{2l}+\mu_{i_{2l-1}}p_{2l-1}+\cdots+\mu_{i_{k}}p_{k}
+\cdots+\mu_{i_{l+\frac{5}{2}}}p_{l+\frac{5}{2}}+\mu_{9}p_{l+\frac{3}{2}}+2\mu_{4}f
-\mu_{6}p_{l+\frac{1}{2}}\nonumber\\
  &&~~~~~~~~~~+\mu_{7}p_{l-\frac{1}{2}}+2\mu_{5}e+\mu_{10}p_{l-\frac{3}{2}}
+\mu_{j_{l-\frac{5}{2}}}p_{l-\frac{5}{2}}+\cdots+\mu_{j_{k'}}p_{k'}+\cdots
+\mu_{j_{1}}p_{1}+\mu_{8}p_{0},\\[6pt]
  \label{der-53}&&D(z)=(2l\mu_{1}+2\mu_{2})z,\\[6pt]
  \label{der-54}&&D(p_{l-\frac{1}{2}})=(l+\frac{1}{2})\mu_{4}p_{l+\frac{1}{2}}
+((l+\frac{1}{2})\mu_{1}+\mu_{2})p_{l-\frac{1}{2}}-(l-\frac{1}{2})\mu_{5}p_{l-\frac{3}{2}}
-(l-\frac{1}{2})!(l+\frac{1}{2})!\mu_{6}z,
  \end{eqnarray*}\begin{eqnarray*}
   \label{der-55}&&D(e)=\mu_{3}p_{2l-1}+\cdots+\frac{k}{2(l-k)}\mu_{i_{k}}p_{k-1}
+\cdots+(-\frac{l+\frac{5}{2}}{5}\mu_{i_{l+\frac{5}{2}}})p_{l+\frac{3}{2}}
-\frac{(l+\frac{3}{2})\mu_{9}}{3}p_{l+\frac{1}{2}}-\mu_{4}h\nonumber\\
  &&~~~~~~~~+(l+\frac{1}{2})\mu_{6}p_{l-\frac{1}{2}}+\mu_{1}e
+(l-\frac{1}{2})\mu_{7}p_{l-\frac{3}{2}}+\frac{(l-\frac{3}{2})}{3}\mu_{10}p_{l-\frac{5}{2}}
+\cdots+\frac{k'}{2(l-k')}\mu_{j_{k'}}p_{k'-1}\nonumber\\
  &&~~~~~~~~+\cdots+\frac{\mu_{j_{1}}}{2(l-1)}p_{0},\\[6pt]
  \label{der-56}&&D(p_{l-\frac{3}{2}})=(l+\frac{3}{2})\mu_{4}p_{l-\frac{1}{2}}
+((l+\frac{3}{2})\mu_{1}+\mu_{2})P_{l-\frac{3}{2}}-(l-\frac{3}{2})\mu_{5}p_{l-\frac{5}{2}}
-\frac{(l+\frac{3}{2})!(l-\frac{3}{2})!}{3}\mu_{9}z,~~~~~\\[6pt]
  \label{der-57}&&D(p_{k'})=(2l-k')\mu_{4}p_{k'+1}+((2l-k')\mu_{1}+\mu_{2})p_{k'}
-k'\mu_{5}p_{k'-1}+\frac{(-1)^{3l-k'+\frac{1}{2}}(2l-k')!k'!}{2(l-k')}\mu_{i_{2l-k'}}z,\\[6pt]
  \label{der-58}&&D(p_{0})=2l\mu_{4}p_{1}+(2l\mu_{1}+\mu_{2})p_{0}
+(-1)^{3l-\frac{1}{2}}(2l-1)!\mu_{3}z,\end{eqnarray*}
  Now we denote by $\lambda_{D}=2l\mu_{1}+2\mu_{2}$ and
  \begin{eqnarray}
\nonumber\!\!\!\!\!\!\!\!
x_{D}\!\!\!&=\!\!\!&\frac{(-\mu_{3})}{2l}p_{2l}+\frac{\mu_{i_{2l-1}}}{2l-2}p_{2l-1}
+\cdots+\frac{-\mu_{i_{k}}}{2(l-k)}p_{k}+\cdots+\frac{\mu_{i_{l+\frac{5}{2}}}}{5}p_{l+\frac{5}{2}}
+\frac{\mu_{9}}{3}p_{l+\frac{3}{2}} \\
\nonumber
\!\!\!\!\!\!\!\!
&+&\mu_{4}f-\mu_{6}p_{l+\frac{1}{2}}-\mu_{7}p_{l-\frac{1}{2}}
-\mu_{5}e+\frac{-\mu_{10}}{3}p_{l-\frac{3}{2}}+\frac{-\mu_{j_{l-\frac{5}{2}}}}{5}p_{l-\frac{5}{2}}
+\cdots+\frac{-\mu_{j_{k'}}}{2(l-k')}p_{k'}+\cdots \\
\nonumber
\!\!\!\!\!\!\!\!
&+&\frac{-\mu_{j_{1}}}{2l-2}p_{1}+\frac{-\mu_{8}}{2l}p_{0}
+\frac{\mu_{1}}{2}h+t_{D}z
\end{eqnarray}
for some $t_{D}\in \mathbb{C}$ associated with $D$. Then it follows that $D(f)=adx_{D}f$,
$D(h)=adx_{D}h$, $D(e)=adx_{D}e$,
$D(z)=adx_{D}z+\lambda_{D}z$, $D(p_{k})=adx_{D}p_{k}+\frac{\lambda_{D}}{2}p_{k}, k=0,1,2,\ldots,2l$.
Let $\delta$ be the linear map from $\mathfrak{g}$ into itself given by (\ref{out-01}), then we
have $D(y)=adx_{D}(y)+\lambda_{D}\delta(y)$ for all $y\in \mathfrak{g}$. The proof is completed.
\hfill$\Box$\vskip7pt

\begin{lemm}\label{xingshi}
Let $\D$ be a 2-local derivation on the conformal Galilei algebra. Then
for every $x,y\in\gi$, there exists a derivation $\D_{x,y}$ of $\gi$ for
which $\D_{x,y}(x)=\D(x)$ and $\D_{x,y}(y)=\D(y)$ and it can be written
as
$$\D_{x,y}=\ad (\sum_{i=0}^{2l}a_i(x,y)p_i
+b_1(x,y)e+b_2(x,y)f+b_3(x,y)h+b_4(x,y)z)+\l(x,y)\d.$$
where $a_i(i=0,1,\dots,2l)$, $b_j(j=1,2,3,4)$, $\l$ are complex-valued
functions on $\gi\times\gi$ and $\d$ is given by Lemma \ref{daozi}.
\end{lemm}

\cl{\bf\S3.\ 2-local derivations on the conformal Galilei algebra}
\setcounter{section}{3}\setcounter{theo}{0}\setcounter{equation}{0}
In this section, we will determine all 2-local derivations on the conformal Galilei algebra.

\begin{lemm}\label{efhz}
Let $\D$ is a 2-local derivation on $\gi$. For any but fixed $x\in\gi$.

(1) If $\D(h)=0$, then $\D_{h,x}=\ad(b_3(h,x)h+b_4(h,x)z)+\l(h,x)\d$;

(2) If $\D(e)=0$, then
    $\D_{e,x}=\ad(a_0(e,x)p_0+b_1(e,x)e+b_4(e,x)z)+\l(e,x)\d$;

(3) If $\D(f)=0$, then
    $\D_{f,x}=\ad(a_{2l}(f,x)p_{2l}+b_2(f,x)f+b_4(f,x)z)+\l(f,x)\d$;

(4) If $\D(z)=0$, then
    $\D_{z,x}=\ad(\sum_{i=0}^{2l}a_i(z,x)p_i
    +b_1(z,x)e+b_2(z,x)f+b_3(z,x)h+b_4(z,x)z)$.
\end{lemm}
\noindent{\it Proof.~}
By Lemma \ref{xingshi}, for any $x\in\{e,f,h,z\}$ we have
$$\D_{x,y}=\ad (\sum_{i=0}^{2l}a_i(x,y)p_i
+b_1(x,y)e+b_2(x,y)f+b_3(x,y)h+b_4(x,y)z)+\l(x,y)\d, $$
where $a_i(i=0,1,\dots,2l)$, $b_j(j=1,2,3,4)$, $\l$ are complex-valued
functions on $\gi\times\gi$ and $\d$ is given by Lemma \ref{daozi}.

When $\D(h)=0$, we have
\begin{align*}
\D(h)&=\D_{h,x}(h)\\
&=[\sum_{i=0}^{2l}a_i(h,x)p_i+b_1(h,x)e
+b_2(h,x)f+b_3(h,x)h+b_4(h,x)z,h]+\l(h,x)\d(h)\\
&=-\sum_{i=0}^{2l}2(l-i)a_i(h,x)p_{i}-b_1(h,x)e+2b_2(h,x)f=0,
\end{align*}
which means $a_i(h,x)=b_1(h,x)=b_2(h,x)=0$ for $i=1,2,\dots,2l$, the proof
 of (1) is completed.

When $\D(e)=0$, we have
\begin{align*}
  \D(e)&=\D_{e,x}(e)\\
&=[\sum_{i=0}^{2l}a_i(e,x)p_i+b_1(e,x)e
+b_2(e,x)f+b_3(e,x)h+b_4(e,x)z,e]+\l(e,x)\d(e)\\
&=-\sum_{i=0}^{2l}ia_i(e,x)p_{i-1}-b_2(e,x)h+2b_3(e,x)e=0,
\end{align*}
which means $a_i(e,x)=b_2(e,x)=b_3(e,x)=0$ for $i=1,2,\dots,2l$, the proof
 of (2) is completed.
Similarly, we can obtain (3).

When $\D(z)=0$, we have
\begin{align*}
\D(z)&=\D_{z,x}(z)\\
&=[\sum_{i=0}^{2l}a_i(z,x)p_i+b_1(z,x)e
+b_2(z,x)f+b_3(z,x)h+b_4(z,x)z,e]+\l(z,x)\d(z)\\
&=\l(z,x)z=0,
\end{align*}
thus we obtain $\l(z,x)=0$,  the proof of (4) is completed.
\hfill$\Box$\vskip7pt

\begin{lemm}\label{sl2}
If $\D$ is a 2-local derivation on $\gi$ such that $\D(h)=\D(e)=0$, then
$\D(f)=0$.
\end{lemm}
\noindent{\it Proof.~}
Since $\D(h)=\D(e)=0$, we assume that
\begin{align}
&\D_{h,x}=\ad(b_3(h,x)h+b_4(h,x)z)+\l(h,x)\d,\label{h}\\
&\D_{e,x}=\ad(a_0(e,x)p_0+b_1(e,x)e+b_4(e,x)z)+\l(e,x)\d,\label{e}
\end{align}
for all $x\in\gi$, where $a_0,b_1,b_3,b_4$ are complex-valued
functions on $\gi\times\gi$.
Take $x=f$ in (\ref{h}) and (\ref{e}) respectively we get
$$\D(f)=\D_{h,f}(f)=[(b_3(h,f)h+b_4(h,f)z,f]+\l(h,f)\d(f)=-2b_3(h,f)f,$$
and
$$\D(f)=\D_{e,f}(f)=[a_0(e,f)p_0+b_1(e,f)e+b_4(e,f)z,f]+\l(e,f)\d(f)
=-2la_0(e,f)p_1+b_1(e,f)h.$$
Comparing the coefficients of above equations, we get
$b_3(h,f)=a_0(e,f)=b_1(e,f)=0$. It concludes that $\D(f)=0$.
\hfill$\Box$\vskip7pt

\begin{lemm}\label{x}
Let $\D$ is a 2-local derivation on $\gi$ such that $\D(h)=\D(e)=\D(f)=0$,
 then for any $x=\sum_{i=0}^{2l}\a_ip_i+k_1e+k_2f+k_3h+k_4z\in\gi$, we have
 $$\D(x)=\l_x(\sum_{i=0}^{2l}\a_ip_i+k_4z), $$
 where $\l_x$ is a complex number depending on $x$.
\end{lemm}
\noindent{\it Proof.~}
Suppose $x=\sum_{i=0}^{2l}\a_ip_i+k_1e+k_2f+k_3h+k_4z\in\gi$.  Using
$\D(h)=\D(e)=\D(f)=0$ and Lemma \ref{sl2}, we have
\begin{align*}
\D(x)&=\D_{h,x}(x)\\
&=[b_3(h,x)h+b_4(h,x)z,x]+\l(h,x)\d(x)\\
&=b_3(h,x)(\sum_{i=0}^{2l}2(l-i)\a_ip_i+2k_1e-2k_2f)
   +\l(h,x)(\frac{1}{2}\sum_{i=0}^{2l}\a_ip_i+k_4z),
\end{align*}
\begin{align*}
\D(x)&=\D_{e,x}(x)\\
&=[a_0(e,x)p_0+b_1(e,x)e+b_4(e,x)z,x]+\l(e,x)\d(x)\\
&=a_0(e,x)((-1)^{l+\frac{1}{2}}(2l)!\a_{2l}z-2lk_2P_1+2lk_3P_0)
  +b_1(e,x)(\sum_{i=0}^{2l}i\a_ip_{i-1}+k_2h-k_3e)\\
  &\ \ \ \ +\l(e,x)(\frac{1}{2}\sum_{i=0}^{2l}\a_ip_i+k_4z),
\end{align*}
and
\begin{align*}
\D(x)&=\D_{f,x}(x)\\
&=[a_{2l}(f,x)p_{2l}+b_2(f,x)f+b_4(f,x)z,x]+\l(f,x)\d(x)\\
&=a_{2l}(f,x)((-1)^{3l+\frac{1}{2}}(2l)!\a_{0}z-2lk_1P_{2l-1}-2lk_3P_{2l})
  +b_2(f,x)(\sum_{i=0}^{2l}(2l-i)\a_ip_{i+1}-k_1h+2k_3f)\\
& +\l(f,x)(\frac{1}{2}\sum_{i=0}^{2l}\a_ip_i+k_4z).
\end{align*}
Compare the above equations, by the arbitrary of $x$ we obtain
$$b_3(h,x)=a_0(e,x)=b_1(e,x)=a_{2l}(f,x)=b_2(f,x)=0, $$ and
$$\l(h,x)=\l(e,x)=\l(f,x). $$ Denote $\l_x=\l(h,x)$ is a constant which
dependent on $x$,
then $$\D(x)=\l_x(\sum_{i=0}^{2l}\a_ip_i+k_4z). $$
\hfill$\Box$\vskip7pt

\begin{lemm}
Let $\D$ be a 2-local derivation on $\gi$ such that $\D(h)=\D(e)=\D(z)=0$,
 then for $e+p_0\in\gi$ and $y\in\gi$, we have
 $$\D_{e+p_0,y}=\ad(a_0(e+p,y)p_0+b_1(e+p,y)e+b_4(e+p,y)c). $$
\end{lemm}
\noindent{\it Proof.~}
Since $\D(h)=\D(e)=0$, by Lemma \ref{sl2} and Lemma \ref{x}, we have
$$\D(e+p_0)=\frac{1}{2}\l_{e+p_0}p_0.$$
Using $\D(z)=0$ and Lemma \ref{efhz}, we obtain
\begin{align*}
\D(e+p_0)=&\D_{z,e+p_0}(e+p_0)\\
=&[\sum_{i=0}^{2l}a_i(z,e+p_0)p_i+b_1(z,e+p_0)e+b_2(z,e+p_0)f
\\&+b_3(z,e+p_0)h+b_4(z,e+p_0)z,e+p_0]\\
=&-\sum_{i=0}^{2l}ia_i(z,e+p_0)p_{i-1}
+(-1)^{3l+\frac{1}{2}}(2l)!\a_{2l}a_i(z,e+p_0)z\\
&+b_2(z,e+p_0)(-h+2lp_1)
+b_3(z,e+p_0)(2e-2lp_0).
\end{align*}
Comparing the coefficients of above equations, we get
$$a_i(z,e+p_0)=0, i=1,\dots,2l, b_2(z,e+p_0)=b_3(z,e+p_0)=0. $$
Then we have $\l_{e+p_0}=0$, and hence $\D(e+p_0)=0$.
For any $y\in\gi$, by Lemma \ref{xingshi} we assume that
\begin{align*}
\D_{e+p_0,y}=&\ad(\sum_{i=0}^{2l}a_i(e+p_0,y)p_i+b_1(e+p_0,y)e
+b_2(e+p_0,y)f\\
&+b_3(e+p_0,y)h+b_4(e+p_0,y)z)+\l(e+p_0,y)\d.
\end{align*}
Then we have
\begin{align*}
\D(e+p_0)&=\D_{e+p_0,y}(e+p_0)\\
&=[\sum_{i=0}^{2l}a_i(e+p_0,y)p_i+b_1(e+p_0,y)e+b_2(e+p_0,y)f\\
&+ b_3(e+p_0,y)h
+b_4(e+p_0,y)z,e+p_0]+\l(e+p_0,y)\d(e+p_0)\\
&=\sum_{i=0}^{2l}-ia_i(e+p_0,y)p_{i-1}+(-1)^{3l+1}(2l)!a_{2l}(e+p_0,y)z
+b_2(e+p_0,y)(-h+2lp_{1})\\
 &+b_3(e+p_0,y)(2e-2(l-k)p_0)+\frac{1}{2}\l(e+p_0,y)p_0=0.
\end{align*}
Thus $\l(e+p_0,y)=b_3(e+p_0,y)=b_2(e+p_0,y)=a_i(e+p_0,y)=0$ for
 $i=1,2,\dots,2l.$
 Therefore
 $$\D_{e+p_0,y}=\ad(a_0(e+p_0,y)p_0+b_1(e+p_0,y)e+b_4(e+p_0,y)z).$$
\hfill$\Box$\vskip7pt

\begin{lemm}\label{zong}
Let $\D$ be a 2-local derivation on $\gi$ such that
$\D(h)=\D(e)=\D(z)=0$,
 then $\D(x)=0$ for any $x\in\gi$.
\end{lemm}
\noindent{\it Proof.~}
Suppose $x=\sum_{i=0}^{2l}\a_ip_i+k_1e+k_2f+k_3h+k_4z\in\gi$, where
$\a_i, k_1, k_2, k_3, k_4\in\C$, where $i=0, 1, \dots, 2l$.
Since $\D(e)=\D(h)=0$, we obtain $\D(f)=0$.
Thus
$$\D(x)=\l_x(\frac{1}{2}\sum_{i=0}^{2l}\a_ip_i+k_4z), $$
for $\l_x\in\C$.
Moreover $\D(c)=0$, we have
$$\D_{e+p_0,y}=\ad(a_0(e+p_0,y)p_0+b_1(e+p_0,y)e+b_4(e+p_0,y)z).$$
Therefore
\begin{align*}
\D(x)=&\D_{e+p_0,x}(x)\\
=&[a_0(e+p_0,x)p_0+b_1(e+p_0,x)e+b_4(e+p_0,x)z,x]\\
=&a_0(e+p_0,x)((-1)^{l+\frac{1}{2}}(2l)!\a_0z-2lk_2p_1-2lk_3p_0)\\
&+b_1(e+p_0,x)(\sum_{i=0}^{2l}i\a_ip_{i-1}+k_2h-k_3e).
\end{align*}

If $k_1=k_2=k_3=0$, i.e. $x=\sum_{i=0}^{2l}\a_ip_i+k_4z$, then we have
\begin{align*}
\D(x)=&\l_x(\frac{1}{2}\sum_{i=0}^{2l}\a_ip_i+k_4z)\\
=&(-1)^{l+\frac{1}{2}}(2l)!a_0(e+p_0,x)\a_0z
+b_1(e+p_0,x)(\sum_{i=0}^{2l}i\a_ip_{i-1}).
\end{align*}
This implies
$$\l_x\a_{2l}=0, $$
$$\frac{1}{2}\l_x\a_{i-1}=b_1(e+p_0,x)i\a_i, $$
and
$$\l_xk_4=(-1)^{l+\frac{1}{2}}(2l)!a_0(e+p_0,x). $$
Let $\a_i$, $k_4$ run
 all integers, we conclude that
 $$\l_x=b_1(e+p_0,x)=a_0(e+p_0,x)=0. $$
 Hence $$\D(x)=0. $$

If $k_4=\a_i=0$ for $i=0,1,\dots,2l$, i.e. $x=k_1e+k_2f+k_3h$,
by Lemma \ref{x}, we have $\D(x)=0$.

If both $k_1,k_2,k_3$ and $k_4, \a_i$ for $i=0,1,\dots,2l$ are not
zero sequences, we assume $k_2\neq0$, then we have
\begin{align*}
b_1(e+p_0,x)&=0, \\
\l_xk_4&=(-1)^{l+\frac{1}{2}}(2l)!a_0(e+p_0,x), \\
\frac{1}{2}\l_x\a_0&=-2lk_3a_0(e+p_0,x), \\
\frac{1}{2}\l_x\a_1&=-2lk_2a_0(e+p_0,x),  \\
\l_x\a_i&=0, i=2,3,\dots,2l.
\end{align*}
Let $\a_0,\a_1$, $k_2,k_3,k_4$ run all integers, we conclude that
$\l_x=a_0(e+p_0,x)=0$, then $\D(x)=0$.
\hfill$\Box$\vskip7pt

Now we give the main theorem concerning 2-local derivations on the conformal Galilei algebra $\gi$.
\begin{theo}
Every 2-local derivation on $\gi$ is a derivation.
\end{theo}
\noindent{\it Proof.~}
Let $\D$ is a 2-local derivation on $\gi$. There exists a derivation $\D_{h,e}$ such that $\D(h)=\D_{h,e}(h)$ and $\D(e)=\D_{h,e}(e)$.
Denote $\D_1=\D-\D_{h,e}$. Then $\D_1$ is a 2-local derivations for
which $\D_1(h)=\D_1(e)=0$. By Lemma \ref{sl2}, we have $\D_1(f)=0$.
By Lemma \ref{x}, we have $\D_1(z)=\l_zz$ for
$\l_z\in\C$. Set $\D_2=\D_1-\l_z\d$. Then $\D_2$ is a 2-local
derivation such that
\begin{align*}
&\D_2(h)=\D_1(h)-\l_z\d(h)=0-0=0,\\
&\D_2(e)=\D_1(e)-\l_z\d(e)=0-0=0,\\
&\D_2(z)=\D_1(z)-\l_z\d(z)=\l_z\d(z)-\l_z\d(z)=0.
\end{align*}
By lemma \ref{zong}, we have $\D_2=\D-\D_{e,h}-\l_z\d\equiv0$. Hence
$\D=\D_{e,h}+\l_z\d$ is a derivation.
\hfill$\Box$

 \end{document}